\documentclass{article}

\usepackage{graphicx}

\usepackage{amsmath}

\usepackage{amsfonts}

\usepackage{amssymb}

\newtheorem{theorem}{Theorem}

\newtheorem{lemma}[theorem]{Lemma}

\newtheorem{proposition}[theorem]{Proposition}

\newenvironment{proof}[1][Proof]{\textbf{#1.} }{\ \rule{0.5em}{0.5em}}

\begin{document}

\title{A note on uniformization of Riemann surfaces by Ricci flow}

\author{Xiuxiong Chen, Peng Lu and Gang Tian}

\date{}

\maketitle

In this note, we clarify that the Ricci flow can be used to give
an independent proof of the uniformization theorem of Riemann
surfaces. The key is the simple observation stated in Lemma
\ref{rotsym}.

When metric $g_0$ has positive curvature on $S^2$ or is a metric
on a surface with genus bigger 1, Hamilton proves that the
volume-normalized Ricci flow $g(t)$ with $g(0)=g_0$ converges to a
metric of constant curvature \cite{Ha}, however the proof in \S10
\cite{Ha} of the fact that gradient shrinking soliton on $S^2$ has
positive constant curvature, uses the uniformization theorem.
There is another proof of the fact which uses the Kazdan-Warner
identity (see p. 131 \cite{CK}), but the proof of the
Kazdan-Warner identity uses the uniformization theorem. When $g_0$
does not have positive curvature on $S^2$, Chow proves that $g(t)$
will have positive curvature for large $t$ \cite{Ch1}, in the
proof the entropy formula is proved using the uniformization
theorem. Later different proofs of entropy estimate, which
replaces the entropy formula, are given without using the
uniformization theorem \cite{Ch2}, \cite{CW}. So the only place,
that the uniformization theorem is used for $g(t)$ to converge to
a metric of constant curvature, is to argue that a gradient
shrinking Ricci soliton on $S^2$ has positive constant curvature.

\vskip .2cm
\begin{lemma} \label{rotsym}
Let $(\Sigma^2,g)$ be a two dimensional complete Riemannian
manifold with non trivial Killing vector field $X$. Suppose $X$
vanishes at $O \in \Sigma$, then $(\Sigma^2,g)$ is rotationally
symmetric.
\end{lemma}

\begin{proof}
Let $\Phi_t:\Sigma \rightarrow \Sigma, t \in (-\infty,\infty)$ be
the isometry group generating by $X$,
$\frac{d}{dt}\Phi_t(x)=X(\Phi_t(x))$. By assumption $\Phi_t(O)=O$
for all $t$. Hence the tangent of $\Phi_t$ induces an oriented
linear isometry
\[(\Phi_t)_*: (T_O \Sigma, g(O)) \rightarrow (T_O \Sigma, g(O)).
\]
Since the oriented linear isometry group of $(T_O \Sigma, g(O))$
is $S^1$ and the map $t \rightarrow (\Phi_t)_*$ is a nontrivial
homomorphism, there is a $t_0 >0$ such that $(\Phi_0)_*
=(\Phi_{t_0})_*$. Note that isometry $\Phi_t$ is determined by
$(\Phi_t)_*$ through minimal geodesic starting at $O$, we have
$\Phi_{t_0} = \Phi_0$. We have shown that there is a nontrivial
$S^1$ isometric action on $(\Sigma,g)$. It is clear that if there
is a ray starting at $O$ then $\Sigma$ is topologically
$\mathbb{R}^2$; If there is not any ray starting at $O$ then
$\Sigma$ is topologically $S^2$.
\end{proof}

%%%%%%%%%%%%%%%%%%%%%%%%%%%%%%%%%%%%%%%%%%%%%%%%%%%%%%%%%%%%%

\vskip .2cm

It was known that a complete gradient steady gradient soliton in dimension
two with positive curvature must be a cigar solution (see the
Editors' footnote on p.241-242 \cite{3CY}). Now using the
rotational symmetry from Lemma \ref{rotsym}, we can prove a
similar proposition for gradient shrinking solitons without using the
uniformization theorem.

\begin{proposition} \label{const cur}
If $g$ is a gradient shrinking soliton on closed surface
$\Sigma^2$, then $g$ has positive constant curvature.
\end{proposition}

\begin{proof} 
Recall the gradient shrinking soliton equation
$R_{ij}=c g_{ij}+\nabla_i \nabla_j f$ for $c>0$ and some smooth
function $f$ on $\Sigma$. We choose $c=1$ below. Since $\Sigma$ is
closed, there is a point $O$ such that $\nabla f (O)=0$. Since
$R_{ij}=\frac{1}{2}Rg_{ij}$ in dimension 2, it follows from
soliton equation that $\nabla f$ is a conformal vector field. Let
$J$ be the almost complex structure on $T\Sigma$ defined by
$90^\circ$ counterclockwise rotation. It is well-known that
$J(\nabla f)$ is a Killing vector field (see, for example, the
Editor's footnote on p.241-242 \cite{3CY}). Hence by Lemma
\ref{rotsym} $g$ is rotationally symmetric $g=dr^2+h(r)^2 d
\theta^2, 0 \leq r \leq A <\infty, 0 \leq \theta \leq 2 \pi $.
From the proof of Lemma \ref{rotsym} we may assume $f=f(r)$. Using
Gauss curvature $K_g=-\frac{h''}{h}$ the soliton equation becomes

\[
-\frac{h''}{h}=1+f'' \qquad -\frac{h''}{h}=1+\frac{h'f'}{h}.
\]

Integrating $ f''=\frac{h'f'}{h}$ we get $f'=ah$ for some constant
$a$. Hence $-\frac{h''}{h}=1+ah'$, multiplying $hh'$ and
integrating over $[0,A]$, we get

\[
\left. -\frac{(h')^2}{2} \right\vert_{0}^{A} =\left. \frac{h^2}{2}
\right \vert_0^A +a \int_0^A h (h')^2 dr.
\]

It follows from the metric $dr^2+h(r)^2 d \theta^2$ is smooth on
$\Sigma$ at $r=0$ and $r=A$ that $h(0)=h(A)=0$ and
$h'(0)=-h'(A)=1$. So $a=0$ and $g$ is Einstein metric.
\end{proof}

\vskip .2cm Hamilton proved that any possible nontrivial soliton
on closed surface must be a gradient shrinking soliton (see the
proof in \S 10.1 \cite{Ha}). So there are no nontrivial solitons
on closed surfaces. \vskip .2cm

In higher dimensions, one may have nontrivial
shrinking Ricci solitons. However, X. C and G. T can prove without
using the uniformization theorem that any shrinking K\"ahler-Ricci
soliton with positive bisectional curvature must be
K\"ahler-Einstein \cite{CT}.

\vskip .1cm
P.L thanks Bennett Chow for discussion about the Editor's note in \cite{3CY}.
Authors are supported by NSF research grants.

\vskip .5cm

\noindent Xiuxiong Chen, Dept of Math, University of Wisconsin, Madison, WI 53706

\vskip .2cm

\noindent Peng Lu, Dept of Math, University of Oregon, Eugene, OR 97403

\vskip .2cm

\noindent Gang Tian, Dept of Math, Princeton University, Princeton, NJ 08544

\end{document}